\documentclass[10pt,journal]{{IEEEtran}}
\usepackage{generic}
\usepackage{cite}
\usepackage{amsmath,amssymb,amsfonts}
\usepackage{algorithmic}
\usepackage{graphicx}
\usepackage{textcomp}
\usepackage{mathtools}
\usepackage{algorithm}
\usepackage[table]{xcolor}
\usepackage{makecell}
\usepackage{setspace}
\usepackage{booktabs,chemformula}
\usepackage{blindtext}
\usepackage{subfig}
\usepackage{makecell}
\def\bw{{\mathbf{w}}}
\def\bx{{\mathbf{x}}}
\def\by{{\mathbf{y}}}
\def\bz{{\mathbf{z}}}
\def\bu{{\mathbf{u}}}

\def\MSD{{\mathrm{MSD}}}

\def\BibTeX{{\rm B\kern-.05em{\sc i\kern-.025em b}\kern-.08em
    T\kern-.1667em\lower.7ex\hbox{E}\kern-.125emX}}

\newtheorem{theorem}{Theorem}

\newtheorem{lemma}{Lemma}

\newtheorem{remark}{Remark}
\newtheorem{assumption}{Assumption}
\def\baselinestretch{0.93}
\begin{document}
 \title{Primal-Dual Strategy (PDS) for Composite Optimization Over Directed graphs }
\author{Sajad Zandi and Mehdi Korki, \IEEEmembership{Member, IEEE}
\thanks{S. Zandi is with faculty of Engineering and IT, University of Technology 
Sydney, Ultimo, 2007, Australia ( email: sajad.zandi@student.uts.edu.au).}
\thanks{M. Korki is with School of Science, Computing and Engineering Technologies, Swinburne University of Technology,  Hawthorn, 3122, Australia ( e-mail: mkorki@swin.edu.au).}}

\maketitle

\begin{abstract}
We investigate the distributed multi-agent sharing optimization problem in a directed graph, with a composite objective function consisting of a smooth function plus a convex (possibly non-smooth) function shared by all agents. While adhering to the network connectivity structure, the goal is to minimize the sum of smooth local functions plus a non-smooth function. The proposed Primal-Dual algorithm (PD) is similar to a previous algorithm \cite{b27}, but it has additional benefits. To begin, we investigate the problem in directed graphs, where agents can only communicate in one direction and the combination matrix is not symmetric. Furthermore, the combination matrix is changing over time, and the condition coefficient weights are produced using an adaptive approach. The strong convexity assumption, adaptive coefficient weights, and a new upper bound on step-sizes are used to demonstrate that linear convergence is possible. New upper bounds on step-sizes are derived under the strong convexity assumption and adaptive coefficient weights that are time-varying in the presence of both smooth and non-smooth terms. Simulation results show the efficacy of the proposed algorithm compared to some other algorithms.

\end{abstract}


\section{Introduction}
\label{sec:introduction}
\IEEEPARstart{I}{n} a decentralized approach, $N$ agents, which only exchange information with their immediate neighbors, are connected to form a multi-agent cooperative network for estimating an unknown parameter through in-network processing.
The goal of this work is to modify and utilize the multi-agent sharing optimization problems \cite{b27}, in a distributed approach.
In recent years, there has been a surge of interest in distributed optimization. It is also a solution to minimize the average/sum of $n$ convex functions by a network of $n$ nodes, where a single agent $i$ has only access to a private function $f_{i}$. This problem frequently arises in formation control \cite{b_intro_13} and non-autonomous power control \cite{b_intro_14}.
The coupled multi-agent sharing optimization problem includes two functions \cite{b27}: the sum of smooth local functions plus a convex (possibly non-smooth) function coupling all agents in the network. Agents aim to find the solution $w_{k}^{\star} \in \mathbb{R}^{Q_{k}}$ to the following problem:
 \begin{equation}\label{eq:1}
    \min_{w_{1} ,...., w_{N}}\sum_{k = 1}^{N}J_{k}(w_{k}) + \mathrm{D}(w),
 \end{equation}
where function $J_{k}(w_{k}):\mathbb{R}^{Q_{k}} \rightarrow \mathbb{R}$ is only known to agent $k$, and non-smooth function $\mathrm{D}(w)$ is known by all agents. Note that $\mathrm{D}(w)=g(\sum_{k = 1}^{N}C_{k}w_{k}) : \mathbb{R}^{M} \rightarrow \mathbb{R} \cup (+\infty)$ is a convex (possibly non-smooth) function, where the matrix $C_{k} \in \mathbb{R}^{M \times Q_{k}}$ is known by agent $k$ only. The coupled multi-agent optimization problem of (\ref{eq:1}) represents the sharing formulation, where the agents (with their various variables) are coupled through $\mathrm{D}(w)$. The applications of problem (\ref{eq:1}) include regression over distributed features \cite{b_intro_01}, \cite{b_intro_02}, dictionary learning over distributed models \cite{b_intro_03}, clustering in graphs \cite{b_intro_04}, smart grid control \cite{b_intro_05}, and network utility maximization \cite{b_intro_06}, to name a few.

The authors in \cite{b27} proposed a decentralized algorithm for (\ref{eq:1}) and established its linear convergence to the global solution. They derived their algorithms based on reformulating (\ref{eq:1}) into an equivalent decentralized solution using an equivalent saddle-point problem. In this work, we use a similar technique to develop and derive a novel decentralized algorithm in a directed graph based on adaptive coefficient weights and a new convergence rate. Weight balance is usually not possible in real-world scenarios. As a result, we aim to improve the optimization and learning algorithms that are suitable for and applicable to directed graphs. The main challenge in dealing with directed graphs is the substitution of construction of doubly-stochastic to either row-stochastic or column-stochastic matrix, which are used as weighted adjacency matrices, in directed graphs. The row-stochasticity of the weight matrix guarantees that all agents reach consensus, while the column-stochasticity ensures optimality \cite{b_intro_07}. Most of the methods for optimization over directed graphs are used to combine the average-consensus methods with optimization algorithms designed for undirected graphs. In contrast to undirected graphs, the applications of a directed graph topology\cite{b43} are wider, and it may also result in lower communication costs and simpler topology. Combining average-consensus techniques created for directed graphs with optimization algorithms made for undirected graphs serves as the inspiration for the current methods for optimization over directed graphs. For example, subgradient-push method, described in  \cite{b_intro_08} and extensively investigated in \cite{b_intro_09}, combines push-sum consensus \cite{b_intro_10} and distributed gradient descent (DGD). Directed-Distributed Gradient Descent (D-DGD)\cite{b_intro_07,b_intro_12}, is a linear technique over directed graphs that is based on surplus consensus \cite{b_intro_11} and DGD. However, due to the diminishing step-size, such DGD-based algorithms converge relatively slowly for general convex functions and strongly convex functions.

\subsection{Related Works}
Information exchange within a multi-agent cooperative network (through coupled neighbors) aims to estimate an unknown parameter. The applications of this cooperative learning and estimation are in various fields such as signal processing, optimization, and control \cite{b1}, \cite{b2}. Information exchange can happen in a distributed or central approach. In the central method, data is collected from the whole network and is processed in the fusion center, then it is sent back to the agents. This method is powerful, however, it still has some limitations such as limited bandwidth and vulnerability of the fusion center to failure. Decentralized methods refer to fully distributed solutions. As an application, the distributed solution method enables parallel processing and the distribution of computational load across multiple agents. When the bandwidth around the central agent is low, the performance of network methods degrades. On the other hand, one advantage of decentralized methods is that convergence occurs quickly while bandwidth is limited. Furthermore, in some sensitive applications, privacy and confidentiality are barriers that prevent sending information to fusion centers. On the other hand, in a distributed model, the agents are only allowed to share information with their immediate neighbors. This can be applied in social, transportation, and biological networks \cite{b1,b44}. 

Many studies have made decentralized optimization algorithms in large-scale networks an interesting topic. For instance, consider \cite{b_intro_15,b_intro_16} with one agent at the center of a centralized topology that may fail or violate privacy requirements. The incremental algorithm \cite{b_intro_17,b_intro_20} may be considered as a replacement for a non-centralized topology based on a directed ring network topology. 
Furthermore, to manage a time-varying network, \cite{b_intro_21} proposes a distributed subgradient (DSG) algorithm with slow speed due to a decaying step-size that is required to obtain a consensus and optimal solution.

In this paper, we propose a new convergence rate for a primal-dual algorithm in a directed graph with adaptive coefficient weights.

Moreover,  if $D(\bw)$ in (\ref{eq:1}) is a non-smooth regularizer across all the agents, a proximal decentralized algorithm structure with a fixed point fits with the desired global solution. On the other hand, if $D(\bw) = 0$, then decentralized primal methods can only converge to a biased solution. As a solution, primal methods are required with a decaying step-size to slow down the convergence rate \cite{b27}. As a result, one of the difficulties in this field is the speed of convergence, for which some methods are affected by decreasing step sizes. To accelerate convergence, we proposed a solution that utilizes strong convexity and the Liptschiz-Continuous gradients case. Some methods employed a constant step size that geometrically converges to an error ball centered on the optimal solution\cite{b10}. The authors in \cite{b_intro_22} achieve geometric convergence with the help of symmetric weights. In \cite{b_intro_23}, the authors combine inexact gradient methods and a gradient estimation technique according to dynamic average consensus, where the underlying graphs must be undirected or weight-balanced according to the proposed methods. The convergence rate of some methods can be limited by decreasing step size, depending on whether the objective function is generally convex or strongly convex.

Recently, the authors in \cite{b_intro_24} proposed a fast distributed algorithm, called DEXTRA, where they have chosen an appropriate step-size which resulted in a linear convergence of the algorithm, while the objective functions are strongly convex with Liptschiz-Continuous gradients.  In this method, if we choose a decreasing step-size, the convergence behavior is slow. However, if the step-size is determined by an interval and the objective functions are strongly convex with the Lipschitz-Continuous gradient, the convergence rate geometrically results in a global optimal point \cite{b41}. In ADD-OPT/push-DIGing \cite{b_intro_25}, the improvement of \cite{b_intro_24} provides a geometrical convergence with a sufficiently small step size. In DEXTRA \cite{b_intro_24} and ADD-OPT/push-DIGing \cite{b_intro_25} algorithms, each node needs to know its out-degree in order to build a column-stochastic weights matrix. In \cite{b_intro_26}, the authors eliminate this requirement by using row-stochastic weights and the need for out-degree knowledge, as each agent locally decides the weights assigned to incoming information. Fast methods over directed graphs share one thing in common: they are all based on push-sum type techniques. This results in a nonlinear algorithm since an independent algorithm is required to asymptotically learn either the right or the left eigenvector, which corresponds to the eigenvalue of 1 \cite{b42}. It does, however, result in additional computations and communications among the agents. 

This work focuses on the case where communication between nodes is directed and time-varying. We analyze the convergence rate with a decreasing step-size for the non-smooth function. First, we discuss the problem when $D(\bw)$ is smooth. Decentralized primal methods in \cite{b8}--\cite{b11} can not achieve an exact solution for the convergence rate unless a degraded step-size is applied. In \cite{b12}, to improve the convergence rate, a decentralized method based on the alternating direction method of multipliers (ADMM) technique has been utilized. To simplify the implementation of the convergence rate, a gradient tracking method was applied in \cite{b13}--\cite{b14}, then, dual-domain and accelerated dual gradient descent were utilized in \cite{b15}. In case, $D(\bw)$ is non-smooth, a proximal gradient method is useful to compute the convergence rate, however, the increasing number of inner loops makes it computationally expensive \cite{b16}. If the non-smooth function is replaced by the indicator function of zero, then, a good convergence rate is attained \cite{b17}, \cite{b18}. For the proximal decentralized algorithm, an asymptotic linear convergence exists even though all the functions are piece-wise linear-quadratic (PLQ). Since both large number of iterations and PLQ\cite{b40} are considered, a good convergence rate can be obtained.

\subsection{Contribution}
In this paper, we develop a modified version of proximal-dual (PD) strategy \cite{b27}, \cite{b27_new} over a distributed network, in a directed graph with adaptive coefficient weights. The purpose of this work is to reach the new convergence rates in a directed graph.

The contributions of this paper are summarized as follows:
\begin{itemize}
     \item In contrast to the proposed methods in the literature \cite{b27}, \cite{b27_new}, which use uniform combination coefficients, the proposed PD algorithm employs adaptive combination coefficients with new step sizes.
    \item To have both a good convergence and a good performance, we optimize the adaptive combination coefficients, in a directed graph, to reduce the load of the communications and the energy of each agent.
    \item Consequently, we derive new bounds for step sizes which result in  a lower squared error. (See Theorem 1)
\end{itemize}
Simulation results also show that the proposed algorithm achieves a faster convergence rate.

The rest of the paper is organized as follows. In Section \ref{Sec.II}, we present the problem formulation of the diffusion strategy and proposed primal-dual diffusion algorithm, and saddle points reformulation. We also provide the implementation of proximal adaptive-then-combine (ATC) for dual-diffusion algorithm and diffusion tracking method. In Section \ref{Sec.III}, we present the proposed proximal dual diffusion Method, i.e., PD-dMVC and we derive the new theorem to find the new step-sizes. In Section \ref{Sec.IV}, we provide the simulation results to evaluate the performance of the proposed algorithm. We conclude the paper in Section \ref{Sec.V}.

\subsection{Notation and Symbols}
All over the paper, all vectors are column vectors, with exception of the regression vector denoted by $\bu$, which is a row vector. $I_{M}$ denotes the identity matrix of size $M \times M$ and $\textbf{1}_{N}$ denotes the $N \times 1$ vector with all entries equal
to one. Lowercase letters denote vectors and uppercase letters denote matrices. $X \geq Y$ and $X > Y$ implies that $X - Y$ is positive definite matrix ( $X$ and $Y$ are both symmetric matrices with the same dimension). $col\{ . \}$ denotes a column vector that is structured by stacking the elements on top of each other. $\|x\|_{C}^{2}=x^{T}Cx$, for a vector $x \in \mathbb{R}^{M}$ and a squared matrix $C \geq 0$. The subdifferential $\partial f(x)$ of a function $f$ at $x$ is the set of all subgradients. The proximal operator of a function $f(x)$ with step-size $\mu$ is $\operatorname{prox}_{\mu f}(x)=\underset{u}{\arg \min } f(u)+\frac{1}{2 \mu}\|x-u\|^{2}$.
The conjugate of a function $f$ is defined as $f^{*}(v)=\sup _{x} v^{\top} x-$ $f(x)$. A differentiable function $f$ is $\delta$-smooth and $\nu$-stronglyconvex if $\|\nabla f(x)-\nabla f(y)\| \leq \delta\|x-y\|$ and $(x-y)^{\top}(\nabla f(x)-$ $\nabla f(y)) \geq \nu\|x-y\|^{2}$, respectively, for any $x$ and $y$.

\section{Problem formulation}\label{Sec.II}

\subsection{Problem Reformulation}

\subsubsection{Primal dual decentralized algorithms}

According to some studies, adaptive then combined methods such as primal-dual recursion have larger step-size stability ranges than non-ATC methods. The primal-dual framework was the first primal-dual interpretation of ATC gradient tracking methods[ref-ref]. When the step-size is kept constant, diffusion methods solve optimization problems in the form of approximate problems that converge to biased solutions. To converge to an optimal solution, primal-based methods require the use of decaying step-size; however, decaying step-size slows down the convergence rate. As a result, primal-dual methods are proposed to obtain unbiased problem-solving and converge to exact minimizers while communications are limited. There are two kinds of problems: those with penalties and those without. The benefit of using a penalty term is that it improves the convergence rate of a decentralized algorithm while aggregate cost is conditioned and strongly convex but individuals are not.
An equivalent approach is to reformulate some problems through a network and try to find interested unknown parameters. In fact, while studying a network with a general algorithmic framework as a solution to the decentralized optimization problem, minimizing a sum of local cost functions is an interesting topic. Furthermore, saddle points as a solution structure include primal-descent, dual ascend, and combine equations. In addition, the cooperation between equations is as follows: gradient-descent followed by gradient-ascend, followed by a combination step.

In this section, we consider some saddle-point problems that can reformulate (\ref{eq:1}). A decentralized algorithm cannot be obtained by directly solving the first equivalent saddle point. As a result, more reformulation is required to derive a decentralized solution in a decentralized approach to solve the problem. The network quantities are introduced as follows:
\begin{align}
&\mathcal{W} \triangleq col\{\bw_{1}, \bw_{2}, ...., \bw_{N}\} \in \mathbb{R}^{Q}, Q \triangleq \sum_{k=1}^{N}Q_{k},\\
&\mathcal{J}(\mathcal{W}) \triangleq \sum_{k = 1}^{N}J_{k}(\bw_{k}),\\
&C \triangleq \left [C_{1}, \cdots, C_{K} \right ] \in \mathbb{R}^{M \times Q}.
\end{align}
Then, the problem (\ref{eq:1}) can be written as follows:
\begin{align}\label{eq:2}
\min_{\mathcal{W}} \mathcal{J}(\mathcal{W}) + g(C\mathcal{W}).
\end{align}
\begin{assumption}(Objective Function)\label{Assump1} \cite{b27}\\
Problem (\ref{eq:2}) has a solution $\mathcal{W}^{\star}$, and the function $\mathcal{J}: \mathbb{R}^{Q} \rightarrow \mathbb{R}$ is $\delta$-smooth and $\nu$-strongly-convex with $0<\nu \leq \delta$. Moreover, the function $g: \mathbb{R}^{E} \rightarrow \mathbb{R} \cup\{+\infty\}$ is a proper lower semi-continuous and  a convex function, and there exists $\mathcal{W} \in \mathbb{R}^{Q}$ such that $C\mathcal{W}$ belongs to the relative interior domain of $g$.
\end{assumption}

\begin{remark}
Because the objective function in (\ref{eq:2}) is strongly convex, the global solution $\mathcal{W}^{\star}$ is unique.
\end{remark} 

Problem (\ref{eq:2}) is equivalent to the following saddle-point reformulation under strong-duality and Assumption (\ref{Assump1}).

\begin{align}\label{eq:3}
\min_{\mathcal{W}}\max_{y} \mathcal{J}(\mathcal{W}) + y^{T}C\mathcal{W} - g^{*}(y),
\end{align}
where $g^{*}$ is the conjugate function of $g$ and $y$ is the dual variable. Due to the coupling effect of all agents by dual variable $y$, which requires a fusion center to compute the dual update, solving (\ref{eq:3}) directly does not lead to a decentralized algorithm. As a result, additional reformulations are required. A local copy of dual variable $y$ at agent $k$, i.e., $y_{k}$, is introduced to solve (\ref{eq:3}) in a decentralized manner. The following quantities are introduced as part of the network:
\begin{align}
\mathcal{Y} &\triangleq col\{y_{1}, y_{2}, ...., y_{N}\} \in \mathbb{R}^{MN},
\mathcal{G}^{*}(\mathcal{Y}) &\triangleq \frac{1}{N}\sum_{k = 1}^{N}g^{*}(y_{k}),\nonumber\\
\mathcal{C}_{d} &\triangleq blkdiag{\{C_{k}\}_{k = 1}^{N}} \in \mathbb{R}^{MN \times Q}.
\end{align}
If we define $\mathcal{D} \in \mathbb{R}^{MN \times MN}$ as a symmetric matrix, we can write:
\begin{align}\label{eq:4}
\mathcal{D}\mathcal{Y} &= 0 \Leftrightarrow y_{1} = .... = y_{N}.
\end{align}
The saddle-point type-2 problem is written as:
\begin{align}\label{eq:5}
\min_{\mathcal{W},\mathcal{X}}\max_{\mathcal{Y}} \mathcal{J}(\mathcal{W}) + \mathcal{Y}^{T}\mathcal{C}_{d}\mathcal{W} + \mathcal{Y}^{T}\mathcal{D}\mathcal{X} - \mathcal{G}^{*}(\mathcal{Y}).
\end{align}
The problem formulation of (\ref{eq:5}) can be solved in a decentralized manner by using an optimal-point $(\mathcal{W}^{*}, \mathcal{X}^{*}, \mathcal{Y}^{*})$.
\begin{lemma}(Optimal Point)
     \textit{Subject to assumption .1 and (\ref{eq:4}), a primal-dual pair $(\mathcal{W}^{*}, y^{*})$ and a fixed point $(\mathcal{W}^{*}, \mathcal{X}^{*}, \mathcal{Y}^{*})$  are optimal if, and only if, they satisfy the optimality conditions for saddle-point 1 and saddle-point 2 as below, respectively.}
\end{lemma}
     \begin{align}
     \textbf{Saddle-Point .1:}(\mathcal{W}^{*} , \mathcal{Y}^{*})\nonumber\\
      - C^{T}y^{*} &= \nabla\mathcal{J}(\mathcal{W}^{*}),\label{eq:9}\\
      C\mathcal{W}^{*} &\in \partial g^{*}(y^{*}),\label{eq:10}\\
     \textbf{Saddle-Point .2:}(\mathcal{W}^{*} , \mathcal{X}^{*} , \mathcal{Y}^{*})\nonumber\\
      - \mathcal{C}_{d}^{T}\mathcal{Y}^{*} &= \nabla\mathcal{J}(\mathcal{W}^{*}),\label{eq:6}\\
      \mathcal{D}\mathcal{Y}^{*} &= 0,\label{eq:7}\\
      \mathcal{C}_{d}\mathcal{W}^{*} + \mathcal{D}\mathcal{X}^{*}  &\in \partial \mathcal{G}^{*}(\mathcal{Y}^{*}).\label{eq:8}
     \end{align}
\begin{lemma}
     \textit{if ($\mathcal{W}^{*} , \mathcal{X}^{*} , \mathcal{Y}^{*}$) fulfills the optimality condition (\ref{eq:6}), (\ref{eq:7}), and (\ref{eq:8}), then it holds that $\mathcal{Y}^{*} = \textbf{1}_{K} \otimes y^{*}$ with ($\mathcal{W}^{*} , y^{*}$) satisfying optimality condition (\ref{eq:9}) and (\ref{eq:10}).}
\end{lemma}
\subsubsection{Decentralized Strategy }
Let $\mathcal{Y}_{-1}$ and $\mathcal{W}_{-1}$ take any arbitrary value and $\mathcal{X}_{-1} = 0$.
\begin{align}
\mathcal{W}_{i} &= \mathcal{W}_{i - 1} - \mu_{\omega}\nabla\mathcal{J}(\mathcal{W}_{i-1}) - \mu_{\omega}\mathcal{C}_{d}^{T}\mathcal{Y}_{i-1},\\
\mathcal{Z}_{i} &= \mathcal{Y}_{i-1} +  \mu_{y}\mathcal{C}_{d}\mathcal{W}_{i} + \mathcal{D}\mathcal{X}_{i - 1},\label{eq:11}\\
\mathcal{X}_{i} &= \mathcal{X}_{i-1} - \mathcal{D}\mathcal{Z}_{i},\\
\mathcal{Y}_{i} &= \operatorname{prox}_{\mu_{y}\mathcal{G}^{*}}(\mathcal{\bar{A}}\mathcal{Z}_{i}),\label{eq:12} 
\end{align}
where $\mathcal{\bar{A}} = \bar{A}\otimes I_{M}$ and $\bar{A} \in \mathbb{R}^{N \times N}$ is left-stochastic or right-stochastic combination matrix. 
\subsubsection{Implementation of proximal ATC for dual-diffusion algorithm}
By rewriting (\ref{eq:11})-(\ref{eq:12}) and removing the variable $\mathcal{X}_{i}$, we obtain:
\begin{align}\nonumber
\mathcal{Z}_{i} - \mathcal{Z}_{i-1} &= \mu_{y}\mathcal{C}_{d}(\mathcal{W}_{i} - \mathcal{W}_{i-1}) + \mathcal{Y}_{i-1} - \mathcal{Y}_{i-2} \\\nonumber  &+ \mathcal{D}(\mathcal{X}_{i - 1} - \mathcal{X}_{i - 2})\\ \nonumber
&= \mu_{y}\mathcal{C}_{d}(\mathcal{W}_{i} - \mathcal{W}_{i-1}) + \mathcal{Y}_{i-1} - \mathcal{Y}_{i-2} \\  &+ \mathcal{D}(-\mathcal{D}\mathcal{Z}_{i-1}).\label{eq:13} 
\end{align}
Rearranging (\ref{eq:13}) yields
\begin{align}\nonumber
\mathcal{Z}_{i} &= (I - \mathcal{D}^{2})\mathcal{Z}_{i - 1} + \mathcal{Y}_{i-1} - \mathcal{Y}_{i-2} \\ &+ \mu_{y}\mathcal{C}_{d}(\mathcal{W}_{i} - \mathcal{W}_{i-1}),
\end{align}
\begin{align}\nonumber
\mathcal{W}_{i} &= \mathcal{W}_{i - 1} - \mu_{\bw}\nabla\mathcal{J}(\mathcal{W}_{i-1}) - \mu_{\bw}\mathcal{C}_{d}^{T}\mathcal{Y}_{i-1}\\
\mathcal{Z}_{i} &= (I - \mathcal{D}^{2})\mathcal{Z}_{i - 1} + \mathcal{Y}_{i-1} - \mathcal{Y}_{i-2} + \mu_{y}\mathcal{C}_{d}(\mathcal{W}_{i} - \mathcal{W}_{i-1}),\label{eq:14}\\
\phi_{i} &= \mathcal{\bar{A}}\mathcal{Z}_{i}\label{eq:16},\\\nonumber
\mathcal{Y}_{i} &= \textbf{Prox}_{\mu_{y}\mathcal{G}^{*}}(\phi_{i}).\nonumber 
\end{align}
\subsubsection{Diffusion Tracking Method}
Let $\mathcal{D} = I - \mathcal{\bar{A}}$ and $\mathcal{\bar{A}} = \mathcal{A}$, then, substituting these new variables into (\ref{eq:14}) results in:
\begin{align}\nonumber
\mathcal{Z}_{i} &= \mathcal{A}(2I - \mathcal{A})\mathcal{Z}_{i-1} + \mathcal{Y}_{i-1} - \mathcal{Y}_{i-2} \\ &- \mu_{y}\mathcal{C}_{d}(\mathcal{W}_{i} - \mathcal{W}_{i-1}).\label{eq:15}
\end{align}
Multiplying (\ref{eq:15}) by $\mathcal{\bar{A}} $ and using (\ref{eq:16}), we have
\begin{align}\nonumber
\phi_{i} = \mathcal{\bar{A}}(\phi_{i-1} + \mathcal{Y}_{i-1} - \mathcal{Y}_{i-2} - \mu_{y}\mathcal{C}_{d}(\mathcal{W}_{i} - \mathcal{W}_{i-1})),
\end{align}
 Hence, the final expressions for the diffusion tracking method are:
\begin{align}
\mathcal{W}_{i} &= \mathcal{W}_{i - 1} - \mu_{\bw}\nabla\mathcal{J}(\mathcal{W}_{i-1}) - \mu_{\bw}\mathcal{C}_{d}^{T}\mathcal{Y}_{i-1}\\
\phi_{i} &= \mathcal{\bar{A}}(\phi_{i-1} + \mathcal{Y}_{i-1} - \mathcal{Y}_{i-2} - \mu_{y}\mathcal{C}_{d}(\mathcal{W}_{i} - \mathcal{W}_{i-1}))
\end{align}

\subsubsection{Adaptive Combination Matrix}

We consider the following definitions and assumptions of the combination matrix:

\begin{align}\nonumber
A &= [a_{l,k}] \in \mathbb{R}^{N\times N},\\
\mathcal{A} &= A\otimes I_{M},
\end{align}

where $a_{l,k}$ is zero if there is no edge connection between agents $k$ and $\ell$. Matrix $A$ is assumed to be a primitive symmetric and left/right-stochastic matrix. $\bar{A} = A $, which is also a primitive symmetric left/right-stochastic matrix. $\mathcal{D} = I - \mathcal{A}$. if the eigenvalues of $\mathcal{\bar{A}}$ belong to (-1,1], many choices for $\mathcal{D}$ exist. 
\begin{assumption}(Combination Matrix)
     \textit{ Matrix $\mathcal{D}$ satisfies condition (\ref{eq:4}) and then below condition yield:}
     \begin{align}\label{Assum 2}\nonumber
         I - \mathcal{D} > 0.
     \end{align}
\end{assumption}
\subsubsection{Theorem (Convergence-rate)}
\begin{theorem}
     \textit{Subject to assumption .1 and .2, if $\mathcal{C}_{d}$ is full row rank and the step-sizes $\mu_{\omega}$ and $\mu_{y}$ are strictly positive and satisfy}
     \begin{align}
     \mu_{\omega} &\leq \frac{1}{2\delta}\\
     \mu_{y} &< \frac{\nu}{2\sigma_{max}^{2}(\mathcal{C}_{d})}
     \end{align}
     \textit{it holds for all $i \geq 0$ and some $C_{o} \geq 0$ that:}\\
     \begin{align}
     \|\tilde{\bw}_{i}\|^{2} \leq \gamma^{i}C_{o}
     \end{align}
     \textit{where}\\
     \begin{align}\nonumber
     \gamma = Max&\{ \frac{1 - \mu_{\omega}\nu( 1 - \mu_{\omega}\delta)}{1 - \mu_{y}\mu_{\omega}\sigma_{max}^{2}(\mathcal{C}_{d})}, 1 - \mu_{\omega}\mu_{y}\lambda_{min}(\mathcal{C}_{d}\mathcal{C}_{d}^{T}),\\ \nonumber & \underline{\sigma}^{2}(\mathcal{A})\|\tilde{x}_{i-1}\|^2\} < 1
     \end{align}
\end{theorem}
     See Appendix \ref{app_1}.
     
    \begin{remark} $\lambda_{min}(\mathcal{C}_{d}\mathcal{C}_{d}^{T})$ is the smallest eigenvalue of $\mathcal{C}_{d}\mathcal{C}_{d}^{T}$, $\underline{\sigma}(\mathcal{A})$ denotes the smallest non-zero singular value of $\mathcal{A}$ and $\sigma_{max}(\mathcal{A})$ is the largest singular value of  $\mathcal{A}$.
    \end{remark}
     \begin{remark} According to Theorem .1, PD algorithm has fast convergence for non-smooth g if $C_{k}$ has full row rank.
     \end{remark}
\subsubsection{ Data Model }
     \begin{enumerate}
     
     \item If we assume that cost function $J(\bw)$ is $\delta$-smooth ,then, following inequality is established,
     \begin{align}
      &\|\nabla J(\bw(i - 1)) - \nabla J(w^{*})\|^2\nonumber \\&\leq \delta\bw^T(i - 1)\Big(\nabla J(\bw(i - 1)) - \nabla J(w^{*})\Big) 
      \end{align}
      \item Strong-convexity is often proposed in adaptation and streaming data to help algorithms against ill-condition by introducing a guard, as well as facilitating conditions to obtain a unique global minimum \cite{b32} - \cite{b33}. If $J(\bw)$ is differentiable, the strong-convexity bound is obtained as
      \begin{align}
      &\tilde{\bw}^T(i - 1)(\nabla J(\bw(i - 1)) - \nabla J(w^{*}))\geq \nu \| \tilde{\bw}(i - 1)\|^2
      \end{align}
     \end{enumerate}
\begin{table*} 	
\caption{Convergence Properties of Distributed Algorithms}
	 \centering {
		\begin{tabular}{lccc}
		\hline
		\rowcolor[gray]{.8}
		 Algorithm &  Support prox.operators &   Rate(convex) &  Step-size(upper bound) \\
		\hline 
 		EXTRA(\cite{b_intro_22}) & Yes  $\displaystyle$ &   $O(\frac{1}{K})$  &  $O\Big(\frac{\upsilon(1 - \rho)}{L^2}\Big)$ \\
 		\rowcolor[gray]{.8}
		Aug-DGM(\cite{b_intro_27})	 & Yes  $\displaystyle$  & converges &  $\min\Big\{\frac{(1-\rho)^2}{10L\rho\sqrt{n}\sqrt{\kappa}},\frac{1}{2L}\Big\}$  \\
		DIGing(\cite{b14})   & Yes  $\displaystyle$  &  $-$  &  $O\Big(\frac{(1-\rho)^2}{\upsilon\kappa^{1.5}\sqrt{n}}\Big)$\\
		\rowcolor[gray]{.8}
		NIDS(\cite{b_intro_28})	 &  Yes $\displaystyle $ &  $O(\frac{1}{K})$ &  $\frac{2}{L}$  \\
		Exact Diffusion(\cite{b_intro_29})  &  Yes   $\displaystyle$ &  $\frac{\kappa^{2}}{1-\rho}$&  $O(\frac{\upsilon}{L^{2}})$ \\
		\rowcolor[gray]{.8}
		Primal-Dual   &  Yes   $\displaystyle$ & $O(\frac{1}{K})$&  $O(\frac{1}{K})$ \\
		Primal-Dual Diffusion  &  Yes  $\displaystyle$ &$O(\frac{1}{K})$& $O(\frac{1}{K})$ \\
		\hline
		\end{tabular}
		}
	\label{tab-1}
\end{table*}

\section{Proximal Decentralized Algorithm}\label{Sec.III}
To manage the non-smooth term $\mathrm{D}(\bw)$, it is defined as
\begin{align}\label{eq:18}
\mathrm{D}(\bw) \triangleq \frac{1}{N}\sum_{k = 1}^{N}\mathrm{D}(\bw_{k})
\end{align}
\subsection{Proximal Dual Diffusion Method}
Following recursion suggested based on the (\ref{eq:17}), for $i = 0,1, ....$, and each agent $k$, if $\mathcal{\bar{A}} = \mathcal{A}, \mathcal{D} = (I - \mathcal{A})$, then

\begin{align}
\bw_{k,i} &= \bw_{k,i - 1} - \mu_{\bw}\nabla J_{k}(\bw_{k,i-1}) - \mu_{\bw}C_{k}^{T}y_{k,i-1}\label{eq:19}\\
\psi_{k,i} &= y_{k,i-1} + \mu_{y}C_{k}\bw_{k,i}\label{eq:20}\\
z_{k,i} &= \phi_{k,i-1} + \psi_{k,i} - \psi_{k,i-1}\label{eq:21}\\
\phi_{k,i} &= \sum_{l \in \mathcal{N}_{k}}\bar{a}_{l,k}z_{l,i}\label{eq:22}\\
y_{k,i} &= \textbf{prox}_{\mu_{k}/Kg^{*}}(\phi_{k,i})\label{eq:23}
\end{align}

\begin{algorithm}\label{alg3}

{\fontsize{7.5}{10}\selectfont

\begin{algorithmic}[1]

\caption{Proximal-Dual diffusion with adaptive coefficient weights-Main Algorithm}

\REQUIRE $ C_{k}$, $ \kappa$

Initialize: $\text{Choose}\hspace{1mm} \phi_{k}(-1) = \psi_{k}(-1) = \theta_{k}(-1) = 0 $, for all $k = 1, 2, ...., N$

\FOR{ Repeat $i \geq 0$ }

\STATE \quad set $\bw_{k}(i) = \bw_{k}(i-1) - \mu_{w}\nabla J_{k}(\bw_{k}(i-1)) - \mu_{w}C_{k}^{T}y_{k}(i-1),$

\STATE \quad set $\psi_{k}(i) = y_k(i-1) + \mu_{y}C_{k}\bw_{k}(i),$

\STATE \quad set $ z_{k}(i) = \phi_{k}(i-1) + \psi_{k}(i) - \psi_{k}(i-1),$

\STATE \quad set $ \phi_{k}(i) = \sum_{\ell \in \mathcal{N}_{k}}a_{\ell,k}z_{\ell}(i),$

\STATE \quad set $y_{k}(i) = Prox_{\mu_{y}/kg^{*}}(\phi_{k}(i))$
 
\textbf{Proximal Computation :}

\STATE \quad set $y_{k}(i) = Prox_{\mu_{y}/kg^{*}}(\phi_{k}(i))$

\STATE \quad set $ \varpi = |\phi_{k}(i)| - \mu_{y}/k$

\STATE \quad set $ \varpi^{'} = max(\varpi, \kappa)exp(-\varpi)$

\IF {$\|\phi_{k}(i)\|_{2} < \varpi^{'}$}

\STATE \quad set $\theta_{1} = F_{k}(i)\sqrt{\|\epsilon_{k}(i)\|_{2}},$

\STATE \quad set $\theta_{2} = \|\phi_{k}(i)\|_{2},$

\STATE \quad set $\theta_{k}(i) = \frac{\theta_{1}}{1 + \theta_{2}}sign(\phi_{k}(i))$

\ELSE

\STATE \quad set $\theta_{k}(i) = \varpi^{'}sign(\phi_{k}(i))$

\ENDIF

\ENDFOR

\end{algorithmic}
}
\end{algorithm}

\begin{table*}
\caption{Performance comparison of algorithms for directed and undirected graphs}
	 \centering{
		\begin{tabular}{lcc}
		\hline
	\rowcolor[gray]{.8}
Algorithm &  Total time (second) &  Number of iterations  \\
		\hline 
 EXTRA(\cite{b_intro_22})&  $19.80$  $\displaystyle$&  $500$ \\
 \rowcolor[gray]{.8}
 Aug-DGM(\cite{b_intro_27})&  $18.58$  $\displaystyle$&  $500$ \\
 DIGing(\cite{b14})&  $18.45$  $\displaystyle$&  $500$ \\
 \rowcolor[gray]{.8}
 NIDS(\cite{b_intro_28})&  $18.31$ $\displaystyle $&  $500$ \\
 Exact Diffusion(\cite{b_intro_29})& $25.62$   $\displaystyle$&  $500$ \\
 \rowcolor[gray]{.8}
 Primal-Dual &  $18.01$  $\displaystyle$ &  $500$\\
 Primal-Dual Diffusion & $17.15$  $\displaystyle$ &  $500$ \\
		\hline
		\end{tabular}
		}
	\label{tab-2}
\end{table*}

The desired variance can be estimated iteratively, and the factor parameter estimated by agent $\ell$ by running the following smoothing filter:
\begin{align}
\gamma_{\ell,k}(i) = (1 - \zeta)\gamma_{\ell,k}(i - 1) + \zeta\|\bw_{\ell}(i-1)\|^{2},
\end{align}
where \{$\bw_{\ell}(i-1)$\}, which is needed to run the recursion at agent k, is computed by agent $k$ at iteration i. In order to remove the necessity of transmitting $\gamma_{\ell,k}(i)$ from agent $\ell$ to its neighbors, the estimation of $\gamma_{\ell,k}(i)$ designed into neighbors of agent $\ell$. It gives the advantage of not sharing the value and overcomes the difficulty of accessing agent k to $\bw_{\ell}(i - 1)$ in the ATC diffusion implementation. This works by the mechanism of receiving from its neighbor $\ell$, we replaced $\bw_{\ell}(i - 1)$ by $\bw_{k}(i - 1)$ since for $i\gg1$. Note that the iterates at various agents approach to an optimal point. With this substitution, agent $k$ can now estimate the variance $\gamma_{\ell}^2(i)$ of its neighbor locally by running a smoothing filter of the following form:

\begin{align}\label{eq:17}
\gamma_{\ell,k}(i) = (1 - \zeta)\gamma_{\ell,k}(i - 1) + \zeta\|\bw_{k}(i-1)\|^{2}.
\end{align}
The above expression provides part of an adaptive construction for the combination weights {$a_{l,k}$}. Besides, adaptive weights are factors to evaluate $\bw_{k}(i)$. 

The proposed Proximal-Dual diffusion algorithm along with its adaptive coefficient weights are described in \textbf{Algorithm 1} and \textbf{Algorithm 2}, respectively.

\begin{algorithm}\label{alg2}

{\fontsize{7.5}{10}\selectfont

\begin{algorithmic}[1]

\caption{Algorithm for Adaptive Coefficient Weights}

\REQUIRE $ \zeta$

Initialize: $\text{Choose}\hspace{1mm} \gamma_{\ell,k}(-1) =  0$, for all $k = 1, 2, ...., N$ and $\ell \in \mathcal{N}_{k}$

\FOR{ Repeat $i \geq 0$ }

\STATE \quad set $\chi_{\ell,k}^{2}(i) = (1 - \zeta)\chi_{\ell,k}^{2}(i-1) + \zeta Max\{ \frac{1 - \mu_{\omega}\nu( 1 - \mu_{\omega}\delta)}{1 - \mu_{y}\mu_{\omega}\sigma_{max}^{2}(\textbf{C}_{d})}, 1 - \mu_{\omega}\mu_{y}\lambda_{min}(\textbf{C}_{d}\textbf{C}_{d}^{T}), \underline{\sigma}^{2}(\mathcal{A})\|\tilde{x}_{i-1}\|^2\} \times \frac{\|\tilde{\bw}_{k}(-1)\|_{I - \mu_{y}\mu_{\omega}C_{d}^T C_{d}}^2 + a_{m}\|\tilde{\by}_{k}(-1)\|_{I - \mu_{y}\mu_{w}\mathcal{C}_{d}\mathcal{C}_{d}^{T}}^2 + a_{m}\|\tilde{\bx}_{k}(-1)\|^2}{1 - \mu_{y}\mu_{w}\sigma^2_(Max)(\mathcal{C}_{d})},\ell\in \mathcal{N}_{k}$

\STATE \quad set $a_{\ell,k}(i) = \frac{exp(\chi_{\ell,k}^{-2}(i))}{\sum_{j \in \mathcal{N}_{k}}exp(\chi_{j,k}^{-2}(i))}, \ell \in \mathcal{N}_{k}$

\ENDFOR

\end{algorithmic}

}
\end{algorithm}

\section{Fast Convergence}\label{Sec.IV}
If $\mathcal{X}_{-1} = 0$, then $\mathcal{X}_{1} = -\mathcal{D}\mathcal{Z}_{1}$ belongs to range space of $\mathcal{D}$. As a result, $\{\mathcal{X}_{i}\}_{i\geq0}$ will always remain in the range space of $\mathcal{D}$. The iterates $(\mathcal{W}_{i}, \mathcal{Y}_{i}, \mathcal{Z}_{i})$ converge to this point $(\mathcal{W}^{*}, \mathcal{Y}^{*}, \mathcal{Z}^{*})$. The error quantities are as follows:
\begin{align}
 \tilde{\mathcal{W}}_{i} &\triangleq \mathcal{W}_{i} - \mathcal{W}^{*} \\
 \tilde{\mathcal{Y}}_{i} &\triangleq \mathcal{Y}_{i} - \mathcal{Y}^{*} \\
 \tilde{\mathcal{Z}}_{i} &\triangleq \mathcal{Z}_{i} - \mathcal{Z}^{*} \\
 \tilde{\mathcal{X}}_{i} &\triangleq \mathcal{X}_{i} - \mathcal{X}^{*}
\end{align}
The expressions in (\ref{eq:19})-(\ref{eq:23}) produce following error recursions:
\begin{align}\nonumber
\mathcal{\tilde{W}}_{i} &= \mathcal{\tilde{W}}_{i - 1} - \mu_{\omega}F_{k}(i)(\nabla\mathcal{J}(\mathcal{W}_{i-1}) - \nabla\mathcal{J}(\mathcal{W}^{*})) \\ &- \mu_{\omega}\mathcal{C}_{d}^{T}\mathcal{\tilde{Y}}_{i-1}\label{eq:50}\\
\mathcal{\tilde{Z}}_{i} &= \mathcal{\tilde{Y}}_{i-1} +  \mu_{y}\mathcal{C}_{d}\mathcal{\tilde{W}}_{i} + \mathcal{D}\mathcal{\tilde{X}}_{i - 1}\label{eq:51}\\
\mathcal{\tilde{X}}_{i} &= \mathcal{\tilde{X}}_{i-1} - \mathcal{D}\mathcal{\tilde{Z}}_{i}\label{eq:52}\\
\mathcal{\tilde{Y}}_{i} &= \textbf{Prox}_{\mu_{y}\mathcal{G}^{*}}(\mathcal{\bar{A}}\mathcal{Z}_{i}) - \textbf{Prox}_{\mu_{y}\mathcal{G}^{*}}(\mathcal{\bar{A}}\mathcal{Z}^{*})\label{eq:53} 
\end{align}

The network mean-square deviation is:\\
  \begin{equation}
  \MSD = \lim_{i \rightarrow \infty} \frac{1}{N}Tr\Big(E\left\{\mathcal{\tilde{Y}}_{i}\mathcal{\tilde{Y}}_{i}^{T}\right\}\Big)
  \end{equation}
  
\section{Convergence properties and performance comparison of algorithms}
Table \ref{tab-1} compares the convergence properties of distributed algorithms in term of supporting prox operators, the rate, and upper bound of the proposed primary-dual diffusion (PDD) and primal-dual(PD) algorithms with those of NIDS, EXTRA and DIGing-ATC, Aug-DGM, Exact diffusion algorithms. As it is seen, the proposed primal-dual diffusion algorithm has the same upper bound of step size as that of the primal-dual algorithm. However, the convergence rate for both PD diffusion and PD is higher than those of the NIDS, EXTRA, DIGing-ATC, Aug-DGM, and Exact diffusion algorithms. Due to the fact shown in Table \ref{tab-1}, PDD and PD require fewer parameters than DIGing, EXTRA and Aug-DGM algorithms, therefore, their upper bounds are simpler. Table \ref{tab-2} compares the performance of the algorithms for directed and undirected graphs in terms of total time and number of iterations. All algorithms were assumed to have the same number of iterations. As shown, the total time per second is different for all algorithms with the same number of iterations. It can be seen that PDD algorithm outperforms the other algorithms.

\section{Simulation Results} \label{Result}

On distributed parameter estimation, simulation for a network with \emph{N} = 20 nodes are presented. Figure 1 (a) depicts the directed network topology used in all simulations. For computing the combination coefficients, the Metropolis rule is used. 
\begin{figure}[!t]
\centering
\subfloat[Network topology]{\includegraphics[height=7cm,width=7cm]{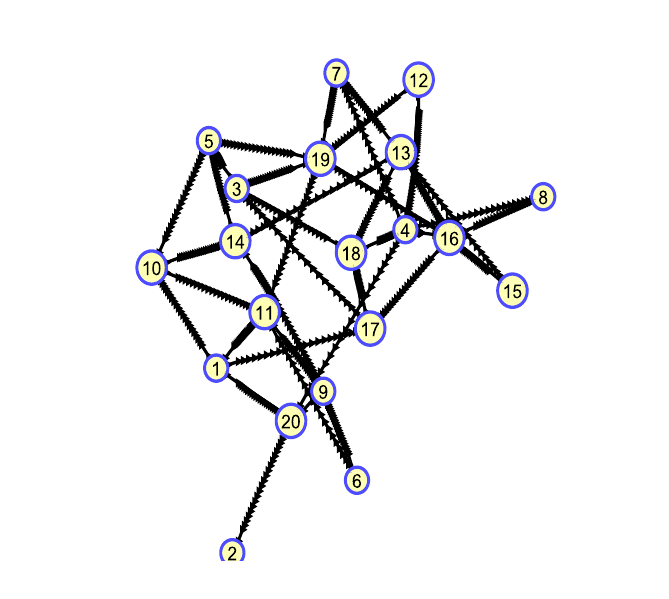}%
\label{fig1:a}}
\hspace{1cm}
\subfloat[Network squared error curve of Primal-dual algorithm without non-smooth function]{\includegraphics[height=7cm,width=8cm]{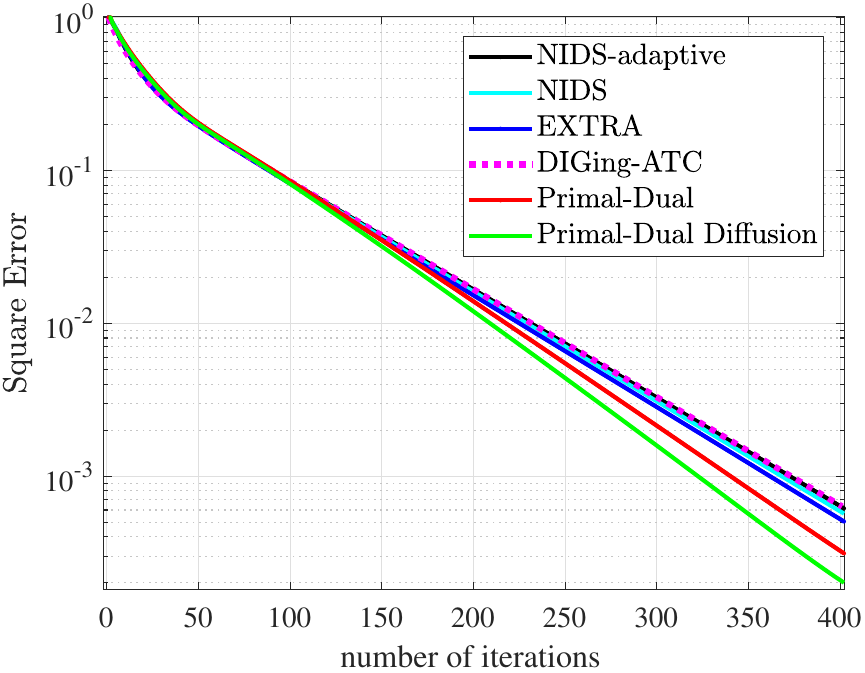}%
\label{fig1:b}}
\caption{Network topology and network squared error curve of Primal-dual algorithm without non-smooth function.}
\label{fig1}
\end{figure}

\subsection{Decentralized compressed sensing} \label{Parameter.Es}

Consider a decentralized compressed sensing problem involving some network agents. Each agent $i \in {1,...., n} $ has its own measurement via $y_{i} = M_{i} x + e_{i} $, where $y_{i} \in \mathbb{R}^{p}$ is an unknown sparse signal and $e_{i} \in \mathbb{R}^{m_{i}} $ is an i.i.d Gaussian noise vector. We consider $\lambda\|x\|$ as a non-smooth function and to simplify the problem, $\lambda$ value is 1. To estimate $x$ as the sparse vector, we use the decentralized algorithm.

\subsubsection{The case without non-smooth function}
Consider the decentralized problem of solving for an unknown signal $x \in \mathbb{R}$ in the absence of a non-smooth function.

\subsubsection{The case with nonsmooth function}
A non-smooth function was considered to solve a decentralized compressed sensing problem.\\
The comparison of the simulation results for the proposed primal-dual diffusion (PDD) and primal-dual (PD) algorithms with those of NIDS-adaptive, NIDS, EXTRA, and DIGing-ATC algorithms are shown Fig. \ref{fig1} (a). Note that in this scenario, we ignored the non-smooth function. It is seen that the proposed PDD outperforms the other algorithms in terms of square error.

In Fig \ref{fig1} (b), we have compared the performance of the proposed PDD algorithm with those of EXTRA and PD
algorithms, where we have a non-smooth function. As it is seen the proposed PDD outperforms EXTRA and PD algorithms.

\begin{figure}[!t]
\centering
\subfloat[without nonsmooth function]{\includegraphics[height=7cm,width=8cm]{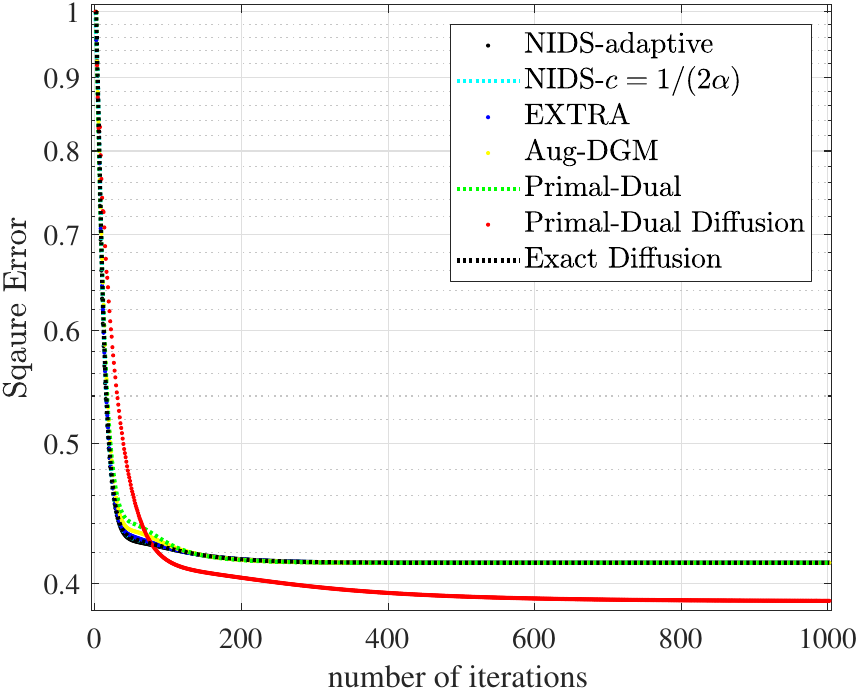}%
\label{fig2:a}}
\hspace{1cm}
\subfloat[with nonsmooth function]{\includegraphics[height=7cm,width=8cm]{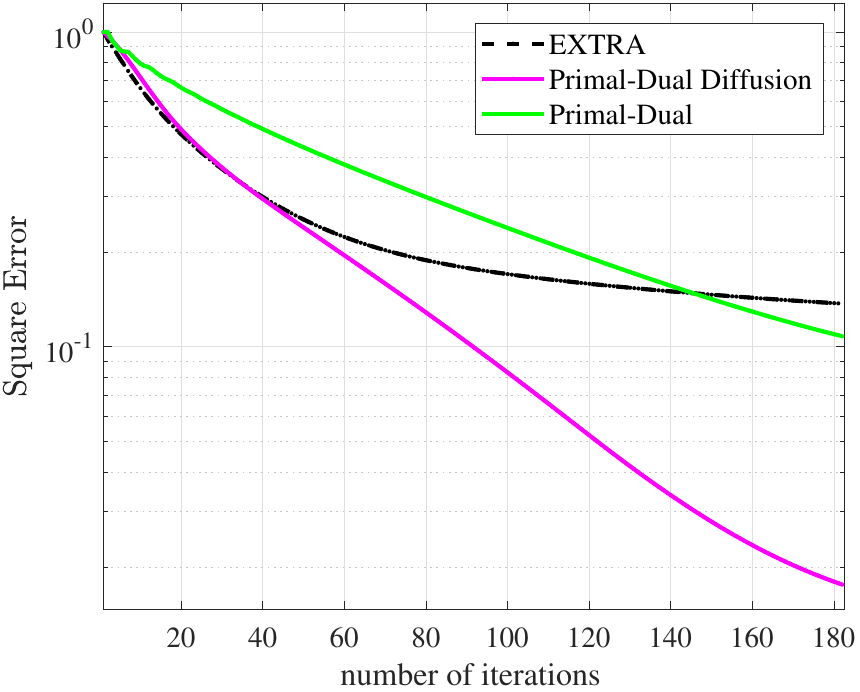}%
\label{fig2:b}}
\caption{Network Squared error curves of Primal-dual algorithm.}
\label{fig2}
\end{figure}

\section{Conclusion}\label{Sec.V}
We have studied the distributed multi-agent sharing optimization problem in a directed graph, with a composite objective function consisting of a smooth function plus a non-smooth function shared through all agents in a network. To solve the problem, some reformulations were applied. A new upper bound on step sizes are obtained, from which a linear convergence can be achieved under the strong convexity assumption. Further, we proposed the adaptive coefficient weights that are time-varying. We show that the proposed algorithm develops the linear convergence to the global minimizer with adaptive coefficient weights and new bound of step sizes.

\appendices
\section{Proof of Theorem .1}\label{app_1}
To solve (\ref{eq:50}, \ref{eq:51}, \ref{eq:52}, \ref{eq:53}) problems, First, we need to define some assumptions, then, square both sides of them to reach optimal step-sizes of PD algorithm and convergence rate equation:
\begin{align}
\tilde{\bw}_{k}(i) &\triangleq \bw_{k}(i) - w^*\\
\tilde{\by}_{k}(i) &\triangleq \by_{k}(i) - y^*\\
\tilde{\bz}_{k}(i) &\triangleq \bz_{k}(i) - z^*\\
\tilde{\bx}_{k}(i) &\triangleq \bx_{k}(i) - x^*
\end{align}
Based on agent formulation, equations are reformulated as follows:
\begin{align}
\bw_{k}(i) &= \bw_{k}(i - 1) - \mu_{\omega}\nabla J(\bw(i - 1)) - \mu_{\omega}\mathcal{C}_{d}^T\by_{k}(i - 1)\\
\tilde{\bw}_{k}(i) &= \tilde{\bw}_{k}(i - 1) - \mu_{\omega}\{\nabla J(\bw(i - 1)) - \nabla J(w^{*})\} \nonumber\\ & \hspace{1em}- \mu_{\omega}\mathcal{C}_{d}^T\tilde{\by}_{k}(i - 1)\label{eq:58}\\
\bz_{k}(i) &= \by_{k}(i - 1) + \mu_{y}\mathcal{C}_{d}\bw_{k}(i) + \mathcal{D}\bx_{i - 1}\\
\tilde{\bz}_{k}(i) &= \tilde{\by}_{k}(i - 1) + \mu_{y}\mathcal{C}_{d}\tilde{\bw}_{k}(i) + \mathcal{D}\tilde{\bx}_{k}(i - 1)\Longrightarrow^{\mathclap{\mathcal{D} = I - A}}\\
\tilde{\bz}_{k}(i) &= \tilde{\by}_{k}(i - 1) + \mu_{y}\mathcal{C}_{d}\tilde{\bw}_{k}(i) + (I - A)\tilde{\bx}_{k}(i - 1)\label{eq:63}\\
\bx_k(i) &= \bx_k(i - 1) - \mathcal{D}\bz_{k}(i)\Longrightarrow^{\mathclap{\mathcal{D} = I - A}}\\
\tilde{\bx}_k(i) &= \tilde{\bx}_k(i - 1) - (I - A)\tilde{\bz}_{k}(i)\label{eq:66}\\
\by_{k}(i) &= Prox(\overline{A}\bz_{k}(i))\\
\tilde{\by}_{k}(i) &= Prox(\overline{A}\bz_{k}(i)) - Prox(\overline{A}\bz^*) \Longrightarrow^{\mathclap{\overline{A} = A}}\\
\tilde{\by}_{k}(i) &= Prox(A\bz_{k}(i)) - Prox(A\bz^*)\label{eq:69}
\end{align}
Now, Squaring of (\ref{eq:58}),(\ref{eq:63}),(\ref{eq:66}), and (\ref{eq:69}),
\begin{align}\nonumber
\|\tilde{\bw}_{k}(i)\|^2 &= \|\tilde{\bw}_{k}(i - 1) - \mu_{\omega}\{\nabla J(\bw_{k}(i - 1)) - \nabla J(w^{*})\}\|^2 \\ & \hspace{1em} -2\mu_{\omega}\tilde{\by}_{k}^T(i - 1)\mathcal{C}_{d}\Big(\tilde{\bw}_{k}(i - 1) \nonumber\\ & \hspace{1em}- \mu_{\omega}\{\nabla J(\bw_{k}(i - 1)) - \nabla J(w^{*})\}\Big)\\
\|\tilde{\bz}_{k}(i)\|^2 &= \|\tilde{\by}_{k}(i - 1) + \mu_{y}\mathcal{C}_{d}\tilde{\bw}_{k}(i)\|^2 + \|(I - A)\tilde{\bx}_k(i - 1)\|^2 \nonumber\\& \hspace{1em} + 2\tilde{\bx}^T_k(i - 1)(I - A)\tilde{\by}_{k}(i - 1) + \mu_{y}\mathcal{C}_{d}\tilde{\bw}_{k}(i)\label{eq:71}\\
\|\tilde{\bx}_{k}(i)\|^2 &= \|\tilde{\bx}_{k}(i - 1)\|^2 + \|(I - A)\tilde{\bz}_k(i)\|^2 \nonumber\\ & \hspace{1em}- 2\tilde{\bx}^T_{k}(i - 1)(I - A)\tilde{\bz}_k(i)\nonumber \\ & =^{\mathclap{(\ref{eq:63})}}\|\tilde{\bx}_{k}(i - 1)\|^2 + \|(I - A)\tilde{\bz}_k(i)\|^2 \nonumber\\ & \hspace{1em}- 2\tilde{\bx}^T_{k}(i - 1)(I - A)(\tilde{\by}_{k}(i - 1)\nonumber\\& \hspace{1em}+ \mu_{y}\mathcal{C}_{d}\tilde{\bw}_{k}(i) + (I - A)\tilde{\bx}_{k}(i - 1))\\
&= \|\tilde{\bx}_{k}(i - 1)\|^2 + \|(I - A)\tilde{\bz}_k(i)\|^2 \nonumber\\ & \hspace{1em}- 2\tilde{\bx}^T_{k}(i - 1)(I - A)(\tilde{\by}_{k}(i - 1)\nonumber\\& \hspace{1em}+ \mu_{y}\mathcal{C}_{d}\tilde{\bw}_{k}(i)) - 2\|(I - A)\tilde{\bx}_{k}(i - 1))\|^2\label{eq:73}
\end{align}
Adding $2\tilde{\bx}^T_{k}(i - 1)(I - A)(\tilde{\by}_{k}(i - 1)+ \mu_{y}\mathcal{C}_{d}\tilde{\bw}_{k}(i))$ and $-\|\bx_{i}\|^2$ to both sides of (\ref{eq:73}) :
\begin{align}\nonumber
  &2\tilde{\bx}^T_{k}(i - 1)(I - A)(\tilde{\by}_{k}(i - 1)+ \mu_{y}\mathcal{C}_{d}\tilde{\bw}_{k}(i)) - \|\bx_{i}\|^2 + \|\bx_{i}\|^2 \\ & = \|\tilde{\bx}_{k}(i - 1)\|^2 + \|(I - A)\tilde{\bz}_k(i)\|^2 \nonumber\\ & \hspace{1em}- 2\tilde{\bx}^T_{k}(i - 1)(I - A)(\tilde{\by}_{k}(i - 1) \nonumber \\ & \hspace{1em}+ \mu_{y}\mathcal{C}_{d}\tilde{\bw}_{k}(i)) - 2\|(I - A)\tilde{\bx}_{k}(i - 1))\|^2 \nonumber\\ & \hspace{1em}+ 2\tilde{\bx}^T_{k}(i - 1)(I - A)(\tilde{\by}_{k}(i - 1) + \mu_{y}\mathcal{C}_{d}\tilde{\bw}_{k}(i)) - \|\bx_{i}\|^2\\
  \nonumber
  \end{align}
  \begin{align}\nonumber
  &2\tilde{\bx}^T_{k}(i - 1)(I - A)(\tilde{\by}_{k}(i - 1)+ \mu_{y}\mathcal{C}_{d}\tilde{\bw}_{k}(i)) \\&= \|\tilde{\bx}_{k}(i - 1)\|^2 + \|(I - A)\tilde{\bz}_k(i)\|^2 \nonumber \\ & \hspace{1em}- 2\|(I - A)\tilde{\bx}_{k}(i - 1))\|^2 - \|\bx_{i}\|^2\label{eq:75}
\end{align}
then, we substitute (\ref{eq:75}) in (\ref{eq:71}), yields,
\begin{align}\nonumber
\|\tilde{\bz}_k(i)\|^2 &= \|\tilde{\by}_{k}(i - 1) + \mu_{y}\mathcal{C}_{d}\tilde{\bw}_{k}(i)\|^2 \nonumber\\ & \hspace{1em}+ \|(I - A)\tilde{\bx}_{k}(i - 1)\|^2+\|\tilde{\bx}_{k}(i)\|^2 \nonumber\\& \hspace{1em}+ \|(I - A)\tilde{\bz}_k(i)\|^2 - 2\|(I - A)(\tilde{\bx}_{k}(i - 1)\|^2\nonumber\\ & \hspace{1em}+ \|\tilde{\bx}_{k}(i))\|^2\nonumber\\ &= \|\tilde{\by}_{k}(i - 1) + \mu_{y}\mathcal{C}_{d}\tilde{\bw}_{k}(i)\|^2 \nonumber\\ & \hspace{1em}+ \|(I - A)\tilde{\bx}_{k}(i - 1)\|^2+\|\tilde{\bx}_{k}(i)\|^2 \nonumber\\& \hspace{1em}+ \|(I - A)\tilde{\bz}_k(i)\|^2 - \|(I - A)(\tilde{\bx}_{k}(i - 1)\|^2\nonumber\\ & \hspace{1em}+ \|\tilde{\bx}_{k}(i))\|^2\label{eq:76}
\end{align}
Both sides of (\ref{eq:76}) should be deducting of $\|(I - A)\tilde{\bz}_k(i)\|^2$, then, yields,
\begin{align}
\|\tilde{\bz}_k(i)\|_{A^2}^2 &= \underbrace{\|\tilde{\by}_{k}(i - 1) + \mu_{y}\mathcal{C}_{d}\tilde{\bw}_{k}(i)\|^2}_{\mathclap{D^{'}}} - \|\tilde{\bx}_{k}(i)\|^2 \nonumber\\ & \hspace{1em}+ \|\tilde{\bx}_{k}(i - 1)\|_{A^2}^2\label{eq:77}
\end{align}
Solving ($D^{'}$):
\begin{align}\nonumber
\|\tilde{\by}_{k}(i - 1) + \mu_{y}\mathcal{C}_{d}\tilde{\bw}_{k}(i)\|^2 &= \|\tilde{\by}_{k}(i - 1)\|^2 + \|\mu_{y}\mathcal{C}_{d}\tilde{\bw}_{k}(i)\|^2 \nonumber\\ & \hspace{1em}+2\tilde{\by}_{k}(i - 1)^T\mu_{y}\mathcal{C}_{d}\tilde{\bw}_{k}(i)\nonumber\\ 
& =^{\mathclap{(\ref{eq:58})}}  \|\tilde{\by}_{k}(i - 1)\|^2 \nonumber\\ & \hspace{1em}+ \|\mu_{y}\mathcal{C}_{d}\tilde{\bw}_{k}(i)\|^2 \nonumber\\ & \hspace{1em}+2\tilde{\by}_{k}(i - 1)^T\mu_{y}\mathcal{C}_{d}\nonumber\\ &\hspace{1em}\times\Big( \tilde{\bw}_{k}(i - 1) \nonumber\\ & \hspace{1em}- \mu_{\omega}\{\nabla J(\bw(i - 1)) \nonumber\\ & \hspace{1em}- \nabla J(w^{*})\} - \mu_{\omega}\mathcal{C}_{d}^T\tilde{\by}_{k}(i - 1)\Big)
\end{align}
similar to (\ref{eq:75}), we can have,
\begin{align}\nonumber
&2\tilde{\by}_{k}^T(i - 1)\mathcal{C}_{d}(\tilde{\bw}_{k}(i - 1) - \mu_{\omega}(\nabla J(\bw_{k}(i - 1)) - \nabla J(w^{*})))\nonumber\\ &=\frac{1}{\mu_{\omega}} (\|\tilde{\bw}_{k}(i - 1) - \mu_{\omega}(\nabla J(\bw_{k}(i - 1)) - \nabla J(w^{*}))\|^2 \nonumber\\ & \hspace{1em}- \|\tilde{\bw}_{k}(i)\|^2 + \|\mu_{\omega}\mathcal{C}_{d}\tilde{\by}_{k}(i - 1)\|^2)\label{eq:79}
\end{align}
if we put (\ref{eq:79}) in ($D^{'}$), yields,
\begin{align}\nonumber
\|\tilde{\by}_{k}(i - 1) + \mu_{y}\mathcal{C}_{d}\tilde{\bw}_{k}(i)\|^2 &= \|\tilde{\by}_{k}(i - 1)\|^2 + \|\mu_{y}\mathcal{C}_{d}\tilde{\bw}_{k}(i)\|^2 \nonumber\\ &\hspace{1em}+ \frac{\mu_{y}}{\mu_{w}}\Big(\tilde{\bw}_{k}(i - 1) \nonumber\\ & \hspace{1em}- \mu_{\omega}(\nabla J(\bw_{k}(i - 1)) - \nabla J(w^{*})) \nonumber\\ &\hspace{1em}- \|\bw_{k}(i)\|^2 + \|\mu_{w}\mathcal{C}_{d}\tilde{\by}_{k}(i - 1)\|^2\Big) \nonumber\\ & \hspace{1em}- 2\mu_{y}\mu_{w}\|\mathcal{C}_{d}\tilde{\by}_{k}(i - 1)\|^2\label{eq:80}
\end{align}
putting (\ref{eq:80}) in (\ref{eq:77}),
\begin{align}
\|\tilde{\bz}_k(i)\|_{A^2}^2 &= \|\tilde{\by}_{k}(i - 1)\|^2 + \|\mu_{y}\mathcal{C}_{d}\tilde{\bw}_{k}(i)\|^2 \nonumber\\ &\hspace{1em}+ \frac{\mu_{y}}{\mu_{w}}\Big( \|\tilde{\bw}_{k}(i - 1) \nonumber\\ &\hspace{1em}- \mu_{\omega}(\nabla J(\bw_{k}(i - 1)) - \nabla J(w^{*}))\|^2 \nonumber\\ &\hspace{1em}- \|\tilde{\bw}_{k}(i)\|^2 + \|\mu_{w}\mathcal{C}_{d}\tilde{\by}_{k}(i - 1)\|^2  \Big ) \nonumber\\ &\hspace{1em}- 2\mu_{y}\mu_{w}\|\mathcal{C}_{d}\tilde{\by}_{k}(i - 1)\|^2 - \|\tilde{\bx}_{k}(i)\|^2 \nonumber\\ &\hspace{1em}+ \|\tilde{\bx}_{k}(i - 1)\|_{A^2}^2\label{eq:81}
\end{align}
Multiplying (\ref{eq:81}) by $a_{m} = \frac{\mu_{w}}{\mu_{y}}$,
\begin{align}
&a_{m}\|\tilde{\bz}_k(i)\|_{A^2}^2 + a_{m}\|\tilde{\bx}_{k}(i)\|^2 = a_{m}\|\tilde{\by}_{k}(i - 1)\|_{I - \mu_{y}\mu_{w}\mathcal{C}_{d}\mathcal{C}_{d}^{T}}^2 \nonumber\\ &\hspace{1em}- \|\tilde{\bw}_{k}(i)\|_{I - \mu_{w}\mu_{y}\mathcal{C}_{d}^{T}\mathcal{C}_{d}}^2 \nonumber\\ &\hspace{1em}+ a_{m}\|\tilde{\bx}_{k}(i - 1)\|^2 \nonumber\\ \hspace{1em}&+ \underbrace{\|\tilde{\bw}_{k}(i - 1) - \mu_{\omega}(\nabla J(\bw_{k}(i - 1)) - \nabla J(w^{*}))\|^2}_{\mathclap{D^{"}}}\label{eq:82}
\end{align}
For solving ($D^{"}$), data model No.6 and data model No.7 needed to be utilize.
\begin{align}
\|\tilde{\bw}_{k}(i - 1) &- \mu_{\omega}(\nabla J(\bw_{k}(i - 1)) - \nabla J(w^{*}))\|^2 \nonumber\\ &= \|\tilde{\bw}_{k}(i - 1)\|^2 + \mu_{\omega}^2\|\nabla J(\bw_{k}(i - 1)) - \nabla J(w^{*})\|^2 \nonumber\\ &\hspace{1em}- 2\mu_{\omega}\tilde{\bw}_{i - 1}^T(\nabla J(\bw_{k}(i - 1)) - \nabla J(w^{*}))
\end{align}
Rely on data model No.6
\begin{align}
\|\tilde{\bw}_{k}(i - 1) &- \mu_{\omega}(\nabla J(\bw_{k}(i - 1)) - \nabla J(w^{*}))\|^2 \nonumber\\ &\leq \|\tilde{\bw}_{k}(i - 1)\|^2 \nonumber\\ &\hspace{1em}+ \mu_{\omega}^2\delta\tilde{\bw}_{i - 1}^T\Big(\nabla J(\bw_{k}(i - 1)) - \nabla J(w^{*})\Big) \nonumber\\ &\hspace{1em}- 2\mu_{\omega}\tilde{\bw}_{i - 1}^T(\nabla J(\bw_{k}(i - 1)) - \nabla J(w^{*}))
\end{align}
and data model No.7,
\begin{align}
\|\tilde{\bw}_{k}(i - 1) &- \mu_{\omega}(\nabla J(\bw_{k}(i - 1)) - \nabla J(w^{*}))\|^2 \nonumber\\ &\leq \|\tilde{\bw}_{k}(i - 1)\|^2 + \mu_{\omega}^2\delta\nu\|\tilde{\bw}_{k}(i - 1)\|^2 \nonumber\\ &\hspace{1em}- 2\mu_{\omega}\nu\|\tilde{\bw}_{k}(i - 1)\|^2
\end{align}
\begin{align}
\|\tilde{\bw}_{k}(i - 1) &- \mu_{\omega}(\nabla J(\bw_{k}(i - 1)) - \nabla J(w^{*}))\|^2 \nonumber\\ &\leq ( 1 -  \mu_{\omega}\nu(2 - \mu_{\omega}\delta))\|\tilde{\bw}_{k}(i - 1)\|^2
\end{align}
then,
\begin{align}
\|&\tilde{\bw}_{k}(i - 1) - \mu_{\omega}(\nabla J(\bw_{k}(i - 1)) - \nabla J(w^{*}))\|^2 \nonumber\\ &\leq ( 1 -  \mu_{\omega}\nu(2 - \mu_{\omega}\delta))\|\tilde{\bw}_{k}(i - 1)\|^2 \nonumber\\
&\leq ( 1 -  \mu_{\omega}\nu(1 - \mu_{\omega}\delta))\|\tilde{\bw}_{k}(i - 1)\|^2 - \mu_{\omega}\nu\|\tilde{\bw}_{k}(i - 1)\|^2\nonumber\\
&\leq ( 1 -  \mu_{\omega}\nu(1 - \mu_{\omega}\delta))\|\tilde{\bw}_{k}(i - 1)\|_{I - \mu_{y}\mu_{\omega}C_{d}^TC_{d}}^2 \nonumber\\
&\leq \frac{ 1 -  \mu_{\omega}\nu(1 - \mu_{\omega}\delta)}{1 - \mu_{y}\mu_{w}\sigma_{Max}^2(C_{d})}\|\tilde{\bw}_{k}(i - 1)\|_{I - \mu_{y}\mu_{\omega}C_{d}^TC_{d}}^2
\end{align}
For positive step-sizes, $a_m > 0$, $1 -  \mu_{\omega}\nu(1 - \mu_{\omega}\delta) > 0$ , then, we have $\mu_{\omega} < \frac{1}{2\delta}$. In addition, $\frac{ 1 -  \mu_{\omega}\nu(1 - \mu_{\omega}\delta)}{1 - \mu_{y}\mu_{w}\sigma_{Max}^2(C_{d})} > 1$, consequently, $\mu_{y} \leq \frac{\nu}{2\sigma_{Max}^2(C_{d})}.$
Finally, (\ref{eq:82}) yields:
\begin{align}
&\|\tilde{\bw}_{k}(i - 1)\|_{I - \mu_{y}\mu_{\omega}C_{d}^TC_{d}}^2 + a_{m}\|\tilde{\bz}_k(i)\|_{A^2}^2 + a_{m}\|\tilde{\bx}_{k}(i)\|^2 \nonumber\\ &= a_{m}\|\tilde{\by}_{k}(i - 1)\|_{I - \mu_{y}\mu_{w}\mathcal{C}_{d}\mathcal{C}_{d}^{T}}^2 + a_{m}\|\tilde{\bx}_{k}(i - 1)\|^2 \nonumber\\ \hspace{1em}&+ \frac{ 1 -  \mu_{\omega}\nu(1 - \mu_{\omega}\delta)}{1 - \mu_{y}\mu_{w}\sigma_{Max}^2(C_{d})}\|\tilde{\bw}_{k}(i - 1)\|_{I - \mu_{y}\mu_{\omega}C_{d}^TC_{d}}^2\label{eq:88}
\end{align}
\begin{assumption}
if $\tilde{\bz}_{k}(i)$ and $\tilde{\bx}_{k}(i - 1)$ lied in the range space of $\mathcal{A}$, then, \begin{align}
\|\tilde{\bx}_{k}(i - 1)\|_{A^2}^2 &\geq  \underline{\sigma}^2(A)\|\tilde{\bx}_{k}(i - 1)\|^2
\end{align}
where $\underline{\sigma}^2(A)$ is the minimum non-zero singular value of $\mathcal{A}$. $\iota_{min}(C_{d}C_{d}^T)$ is assumed as the smallest eigenvalue of $C_{d}C_{d}^T$ which is positive-definite and $\iota_{min}(C_{d}C_{d}^T)I \leq C_{d}C_{d}^T$. 
\begin{align}
    \|\tilde{\by}_{k}(i - 1)\|&_{I - \mu_{y}\mu_{w}\mathcal{C}_{d}\mathcal{C}_{d}^{T}}^2 \nonumber\\&\leq (1 - \mu_{y}\mu_{\omega}\iota_{min}(C_{d}C_{d}^T))\|\tilde{\by}_{k}(i - 1)\|^2
\end{align}
\end{assumption}
\begin{remark}
Squaring both sides of (\ref{eq:69}), introduced below equation,
\begin{align}\nonumber
    \|\tilde{y}_{i}\|^2 &= \|Prox(\mathcal{\bar{A}}\mathcal{Z}_{i}) - Prox(\mathcal{\bar{A}}\mathcal{Z}^{*})\|^2
    \\&\leq\|\mathcal{\bar{A}}\tilde{\mathcal{Z}_{i}}\|^2
     = \|\tilde{\mathcal{Z}_{i}}\|_{\mathcal{\bar{A}}^2}^2\nonumber\\&\leq^{(\mathcal{\bar{A}} = \mathcal{A})}\|\tilde{\mathcal{Z}_{i}}\|_{\mathcal{A}^2}^2
\end{align}
\end{remark}
Based on Assumption 3. and Remark 6., (\ref{eq:88}) changed and yields,
\begin{align}
&\|\tilde{\bw}_{k}(i - 1)\|_{I - \mu_{y}\mu_{\omega}C_{d}^TC_{d}}^2 + a_{m}\|\tilde{\by}_k(i)\|^2 + a_{m}\|\tilde{\bx}_{k}(i)\|^2 \nonumber\\ &= \gamma_1\tilde{\bw}_{k}(i - 1)\|_{I - \mu_{y}\mu_{\omega}C_{d}^TC_{d}}^2 \nonumber\\ \hspace{2em}&+ a_{m}\gamma_2\|\tilde{\by}_{k}(i - 1)\|_{I - \mu_{y}\mu_{w}\mathcal{C}_{d}\mathcal{C}_{d}^{T}}^2 \nonumber\\ \hspace{2em}&+ a_{m}\gamma_3\|\tilde{\bx}_{k}(i - 1)\|^2\end{align} where
\begin{align}
\gamma_1 &\triangleq \frac{ 1 -  \mu_{\omega}\nu(1 - \mu_{\omega}\delta)}{1 - \mu_{y}\mu_{w}\sigma_{Max}^2(C_{d})}\\
\gamma_2 &\triangleq (1 - \mu_{y}\mu_{\omega}\iota_{min}(C_{d}C_{d}^T))\\
\gamma_3 &\triangleq \underline{\sigma}^2(A)
\end{align}

\begin{IEEEbiography}[{\includegraphics[width=1in,height=1in,clip,keepaspectratio]{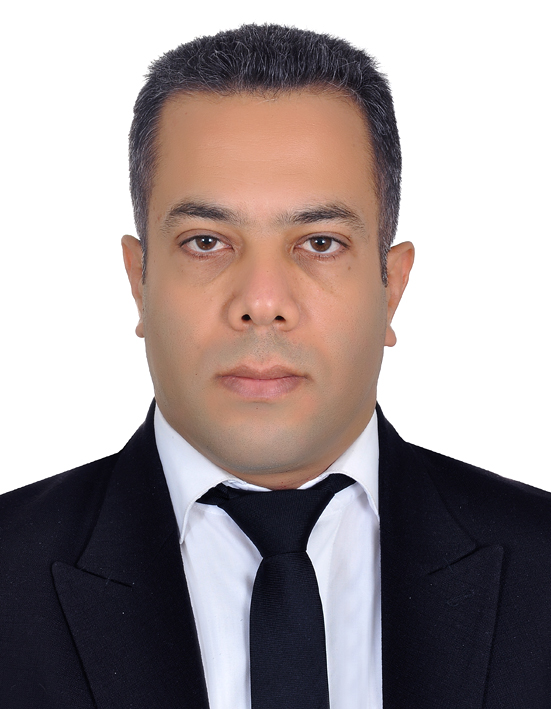}}]{%
Sajad Zandi}
completed his Master of Science degree in Electrical Engineering from Malayer University, Iran, in 2017. Since May 2024, he has been a PhD student at the School of Mechanical and Mechatronics Engineering, University of Technology Sydney (UTS), Australia. His research focuses on adaptive and statistical signal processing, cooperative learning, multi-agent networking, distributed adaptive estimation, distributed active noise control, information-theoretic data, signal processing for communications, and adaptive optimization.
\end{IEEEbiography}

\begin{IEEEbiography}[{\includegraphics[width=1in,height=1.25in,clip,keepaspectratio]{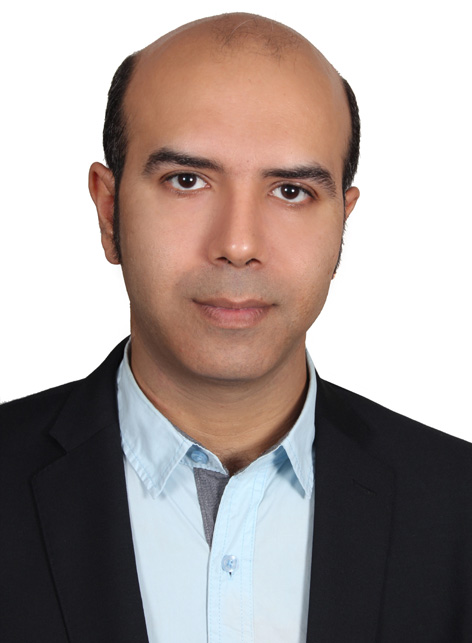}}]{%
Mehdi Korki}
received the bachelor degree in electrical engineering from Shiraz University, Shiraz, Iran, in 2001 and the master and Ph.D. degrees in electrical engineering from Swinburne University of Technology, Melbourne, Australia, in 2012 and 2016, respectively. He is currently a Lecturer at the School of Science, Computing and Engineering Technologies, Swinburne University of Technology, Melbourne, Australia. His research interests include statistical signal processing, data privacy, wireless sensor networks and distributed signal processing.
\end{IEEEbiography}
\end{document}